\documentclass{jande}
\usepackage{jande}
\setVolume{X}{X} 
\setIndex{20XX}{XX}{XXX} 

\usepackage{graphicx}
\usepackage{epstopdf}
\usepackage{cite}

\theoremstyle{definition}  
\theoremstyle{remark}      
\theoremstyle{plain}       

\theoremstyle{remark}

\theoremstyle{definition}

\titlemark{Alexander N. Pchelintsev / Journal of Applied Nonlinear Dynamics 9(2) (2020) 207--221}
\authormark{Alexander N. Pchelintsev / Journal of Applied Nonlinear Dynamics 9(2) (2020) 207--221}

\newcommand{\reff}[1]{(\ref{#1})}
\newcommand{\nl}{\medskip\\}
\newcommand{\disp}{\displaystyle}
\newcommand{\tr}{^{\scriptsize\mbox{T}}}
\newcommand{\indrk}{\{\mbox{\scriptsize rk4}\}}
\newcommand{\dtrk}{\Delta t^{\indrk}}
\newcommand{\eps}{\varepsilon}
\newcommand{\epsrk}{\eps^{\indrk}}

\begin{document}

\title{\first{An Accurate Numerical Method and Algorithm for Constructing Solutions of Chaotic Systems}}

\setcounter{footnote}{1}
\author{\noindent\large Alexander N. Pchelintsev\footnote{Corresponding author.\\Email address:
pchelintsev.an@yandex.ru}}

\address{\normalsize Department of Higher Mathematics, Tambov State Technical University,\\[-2mm]
ul. Sovetskaya 106, Tambov, 392000, Russia}

\abstract{
\begin{table}[h!]
\vspace*{-5mm} \doublerulesep 0.05pt \tabcolsep 7.8mm \vspace*{2mm}
\setlength{\tabcolsep}{7.5pt}
\hspace*{-2.5mm}\begin{tabular*}{171.45mm}{r|||||l}
\multicolumn{2}{l}{\rule[-6pt]{171.45mm}{.01pt}}\\
\parbox[t]{6cm}{\small
\vspace*{.5mm}
\hfill {\bf Submission Info}\par
\vspace*{2mm}
\hfill Communicated by Referees\par
\hfill Received DAY MON YEAR \par
\hfill Accepted DAY MON YEAR\par
\hfill Available online DAY MON YEAR\par
\noindent\rule[-2pt]{6.3cm}{.1pt}\par
\vspace*{2mm}
\hfill {\bf Keywords}\par
\vspace*{2mm}
\hfill Attractor\par
\hfill Power series\par
\hfill Calculation of Lyapunov exponents}
&
\parbox[t]{10.5cm}{
\vspace*{.5mm}
{\normalsize\bf Abstract}\par
\renewcommand{\baselinestretch}{.8}
\normalsize \vspace*{2mm} {\small In various fields of natural science,
the chaotic systems of differential equations are considered more than 50 years.
The correct prediction of the behaviour of solutions of dynamical model equations
is important in understanding of evolution process and reduce uncertainty. However,
often used numerical methods are unable to do it on large time segments.
In this article, the author considers the modern numerical method and algorithm for
constructing solutions of chaotic systems on the example of tumor growth model.
Also a modification of Benettin's algorithm presents for calculation
of Lyapunov exponents.}
\par
\par
\par
}\\
&\\
\multicolumn{2}{l}{\rule[15pt]{171.45mm}{.01pt}}\\\end{tabular*}
\vspace*{-7mm}
\end{table}
}

\maketitle
\thispagestyle{fancy}

\renewcommand{\baselinestretch}{1}
\normalsize

\section{Introduction}
\label{intro}
In 1963, Lorenz considered in \cite{Lorenz} the dynamical system
$$
  \dot{x}_1=\sigma(x_2-x_1),\:\:\dot{x}_2=rx_1-x_2-x_1x_3,\:\:
  \dot{x}_3=x_1x_2-bx_3,
$$
from a model describing a Rayleigh-Benard convection.
At $\sigma=10$, $r=28$ and $b=8/3$ in this system, there is the chaotic
behaviour of solutions, i.e. the solutions are unstable and at the same
time bounded. As it is known from the classical work \cite{Nemytskii},
if the solutions are limited for time $t\ge0$, then the limit set exists.
The trajectories of the dynamical system are attracted to it for
$t\rightarrow\infty$. After approx 13 years after the Lorenz article, the
hypothesis about the structure of the Lorenz attractor was formulated
in \cite{Guckenheimer,AfraimovichBykov,Williams,Kaplan}. It was based on
computational experiments.

Let us show the several dynamic systems with chaotic behaviour of trajectories
at the last 50 years:

\begin{enumerate}
\item In the article \cite{Cook}, the authors described the chaotic change
in time of the magnetic poles of Earth (the Rikitake system).

\item Tyson \cite{Tyson} described the scheme and proposed
a modified equation of Oregonator. It reflects the features of the self-oscillating
chemical reaction of Belousov-Zhabotinsky.

\item Vallis modelled in \cite{Vallis86,Vallis88} the temperature
oscillations in the eastern and western parts of the equatorial region of ocean
which have a strong influence on the global climate of Earth.

\item In \cite{Sprott94,Sprott97}, Sprott showed the nineteen 3rd order
ordinary differential equations (ODEs) and 3rd order ODE which exhibits
chaos. The similar systems were later \cite{Wei,WangChen}.

\item Stenflo \cite{Stenflo} received a system describing the evolution of
amplitude acoustic gravity waves in a rotating atmosphere. The Lorenz--Stenflo
system reduced to the Lorenz system when the parameter tied with the flow rotation
is equal to zero.

\item Chen considered \cite{Chen,Ueta} the controlled Lorenz equation
with a linear feedback controller and received a new system with a double
scroll chaotic attractor.

\item Not so long ago the nonlinear economic systems (e.g. \cite{Magnitskii})
appeared, where there is chaos.

\item The canonical Gause--Lotka--Volterra model for describing of human
sequential memory dynamics is considered in
\cite{AfraimovichGong,AfraimovichVarona}.
\end{enumerate}

The Lorenz system and all of these systems are united not only by the chaotic
behaviour of solutions, but by type of nonlinearities in the right-hand side
of equations. These models have the quadratic nonlinearities. The authors
\cite{Dudkowski} are presented a detailed analysis of the hidden attractors in
some of them.

Many researchers used the classical numerical methods to study the attractors
of dynamical systems. For example, the explicit Euler scheme with
the central-difference scheme \cite{Lorenz}, the Adams method
\cite{Yorke,Yao}, the higher derivatives scheme \cite{Sparrow} and the Runge--Kutta
methods \cite{Cook,Magnitskii,Kaloshin,Sarra}. This methods cannot be used to build
of correct prediction due to the unstability of solutions at a given time segment,
since the global calculating error grows by increasing of time
(the attractors are examined on large time segments). It noted Lorenz in his report
\cite{LorenzButterfly} (the butterfly effect), but such error is limited by the
diameter of sphere, containing an attractor.

Now there are methods that the accumulation of errors is not as great as it was in
the classical methods. Motsa \cite{Motsa,MotsaDlamini} presented a the
piecewise-quasilinearization and multistage spectral relaxation methods which
are based on the Chebyshev spectral method to solve the system and iteration schemes
at each subinterval of integration. In the article \cite{Eftekhari}, the authors used
the differential quadrature method with a similar idea to the solution of system of
ODEs. Another used method is the homotopy-perturbation method \cite{Chowdhury}.

In these methods, the main problems are the choice of integration step and
calculation error of the numerical method.

Starting around the 1960s, the method of power series is starting to
develop for applied computing. Gibbons in \cite{Gibbons} considered the main
types of right-hand sides of ODEs and corresponding computational formulas.
Today this idea was generalized in a recursive procedure (called as automatic
differentiation) to compute the values of the derivatives for power series
\cite{Rall}. An advantage over the general Taylor series method is that
the calculations can be constructed by fast formulas in comparison to
the direct symbolic differentiation of right-hand sides of nonlinear ODEs which
requires a lot of computer memory for high-precision calculations. The method
of power series in \cite{Hashim,Abdulaziz,Sawalha} is applied as the Adomian
decomposition method (ADM). The Clean Numerical Simulation (CNS)
\cite{Liao1,Liao2,Liao3,Liao4,Liao5,Liao6} is based on the method of power series
at arbitrary-order and used the multiple-precision data, plus a check of solution
by means of an additional computation using even smaller numerical noises.

In the FGBFI-method (the firmly grounded backward-forward integration method)
\cite{Pchelintsev,LoziPchelintsev,LozPogPchel}, the authors have taken into
account the above shortcomings of numerical methods used for constructing
solutions of chaotic type, i.e.:
\begin{enumerate}
\item The recurrence relations for calculating of the coefficients of expansion
of local solutions in a power series are received for any dynamic system with
quadratic nonlinearities in the general form.

\item The convergence of the power series is studied. The authors derived a
simple formula of calculating the length of the integration step in
the general form (e.g., it distinguishes the FGBFI-method from CNS).

\item The criteria for checking the accuracy of the approximate chaotic
solution are obtained. There are the control of accuracy and configuration of
obtained approximate solution of a dynamical system with the forward and backward
time which makes the reliability of the numerical method (the degrees of piecewise
polynomials, the value of the maximum step of integration, etc.).
\end{enumerate}

In this paper, the author considers the FGBFI-method for constructing solutions of chaotic
biological system \cite{Mansilla} (the model of tumor growth).
The main advantage of this method is what it allows to produce a more accurate
research of the behaviour of solutions of dynamical systems in very large time segments.
Let us note that the FGBFI-method can be used in the encryption system,
constructed by means of continuous-time chaotic systems \cite{Arroyo},
and also for verification of approximate periodic solutions of continuous nonlinear
dynamical systems \cite{Luo1,Luo2}.

\section{Method of finding of approximate solutions describing the tumor growth}
\label{sec:2}

Let us consider the model developed in the article \cite{Mansilla}:
\begin{eqnarray}
   \left\{\begin{array}{l}
      \dot{x}_1=2Nx_1-x_1^2-Hx_1x_3,\\
      \dot{x}_2=(4-I)x_2+0.5x_1^2-0.14x_2^2-0.5Hx_2x_3+0.001x_3^2,\\
      \dot{x}_3=-Ix_3+0.07x_2^2+0.5Hx_2x_3-0.002x_3^2,
   \end{array}\right.\label{tumorsys}
\end{eqnarray}
where $x_1(t)$, $x_2(t)$ and $x_3(t)$ are a population of proliferating tumor
cells in the avascular, vascular and metastasis phases, respectively;
$N$, $H$ and $I$ are some numbers. The essence of the system parameters:
$N$ is a population of normal cells, $H$ is a population of the host
cells, and $I$ is a population of immune cells (T lymphocytes (CTL)
and natural killer (NK) cells). In this system, there is the
chaotic solutions for certain values of the parameters.

The right side of this system has the quadratic nonlinearities. Then we
can apply the FGBFI-method described in the articles
\cite{Pchelintsev,LoziPchelintsev,LozPogPchel} to construct an accurate
prediction of solutions in a given time segment.

For this purpose, we rewrite system in the vector form \cite{LozPogPchel}
\begin{eqnarray}
   \dot{X}=AX+\Phi(X),\label{dynrsys}
\end{eqnarray}
where
$$
\begin{array}{c}
X(t)=\left[x_1(t)\:x_2(t)\:x_3(t)\right]\tr,\:
\Phi(X)=\left[\varphi_1(X)\:\varphi_2(X)\:\varphi_3(X)\right]\tr,\nl
\varphi_p(X)=\langle Q_p X,X \rangle,\:\:\:p=\overline{1,3},
\:\:\:A=\left[
        \begin{array}{ccc}
           2N &   0 &  0\\
            0 & 4-I &  0\\
            0 &   0 & -I
         \end{array}
  \right],\nl
  Q_1=\left[
         \begin{array}{rrr}
            -1 & 0 & -H\\
             0 & 0 &  0\\
             0 & 0 &  0
         \end{array}
      \right],\:\:
  Q_2=\left[
         \begin{array}{rrr}
            0.5 &     0 &     0\\
              0 & -0.14 & -0.5H\\
              0 &     0 & 0.001
         \end{array}
      \right],\nl
  Q_3=\left[
         \begin{array}{rrr}
            0 &    0 &      0\\
            0 & 0.07 &   0.5H\\
            0 &    0 & -0.002
         \end{array}
      \right].
\end{array}
$$

We expand the solution as
\begin{eqnarray}
   x_1(t)=\sum_{i=0}^\infty\alpha_{1,i} t^i,\:\:\:
   x_2(t)=\sum_{i=0}^\infty\alpha_{2,i} t^i,\:\:\:
   x_3(t)=\sum_{i=0}^\infty\alpha_{3,i} t^i,\label{powser}
\end{eqnarray}
where $x_1(0)=\alpha_{1,0}$, $x_2(0)=\alpha_{2,0}$ and
$x_3(0)=\alpha_{3,0}$ are initial conditions.

The formulas for calculating of the coefficients obtained
as follows: the multiplications of phase coordinates are
assigned by the sums
$$
  \begin{array}{c}
     \disp
     x_1^2=x_1\cdot x_1\:\Rightarrow\:
       r_{1,i}=\sum_{j=0}^{i}\alpha_{1,j}\alpha_{1,i-j},\:\:
     x_2^2\:\Rightarrow\:r_{2,i}=\sum_{j=0}^{i}\alpha_{2,j}\alpha_{2,i-j},\nl
     \disp
     x_3^2\:\Rightarrow\:r_{3,i}=\sum_{j=0}^{i}\alpha_{3,j}\alpha_{3,i-j},\:\:
     x_1x_3\:\Rightarrow\:r_{4,i}=\sum_{j=0}^{i}\alpha_{1,j}\alpha_{3,i-j},\nl
     \disp
     x_2x_3\:\Rightarrow\:r_{5,i}=\sum_{j=0}^{i}\alpha_{2,j}\alpha_{3,i-j}.
  \end{array}
$$
Then the relations for calculating the coefficients of the series
are
\begin{eqnarray}
   \left\{\begin{array}{l}
      \alpha_{1,i+1}=\dfrac{2N\alpha_{1,i}-r_{1,i}-Hr_{4,i}}{i+1},\nl
      \alpha_{2,i+1}=\dfrac{(4-I)\alpha_{2,i}+0.5r_{1,i}-0.14r_{2,i}-
                     0.5Hr_{5,i}+0.001r_{3,i}}{i+1},\nl
      \alpha_{3,i+1}=\dfrac{-I\alpha_{3,i}+0.07r_{2,i}+0.5Hr_{5,i}
                            -0.002r_{3,i}}{i+1}
   \end{array}\right.\label{recrel}
\end{eqnarray}
for $i=0,\:1,\:2,...$ by analogy with \cite{Pchelintsev,LoziPchelintsev}.
This formulas is simpler and faster for calculating than in ADM.

Since the criteria for checking the accuracy of the approximate
chaotic solutions require to go in the backward time repeatedly, then
we need to have a guaranteed estimation of a region of
convergence for given $(\alpha_{1,0},\:\alpha_{2,0},\:\alpha_{3,0})$.
It is usually assumed in some articles (e.g. \cite{Wang})
that the integration step is given and does not change in a calculating
experiment in the nonlinear case, or at all not justified.
We can research the asymptotic behaviour $\alpha_{1,i}$, $\alpha_{2,i}$
and $\alpha_{3,i}$ to determine the integration step, but this question
is poorly investigated today for nonlinear recurrence relations unlike
the linear case \cite{Mezzarobba}.

In the article \cite{LozPogPchel}, the authors proved the theorem
about estimation a region of convergence for the ODEs with
any quadratic nonlinearities. In particular, in this case
(for $I\ge0$, $H>1$, $N>0$)
$$
  \begin{array}{c}
     \|A\|=\|A\|_1=\max\{2N,\:|4-I|,\:I\},\:\:\:
     \|Q_1\|=\|Q_1\|_1=H,\nl
     \|Q_2\|=\|Q_2\|_1=0.5H+0.001,\:\:\:
     \|Q_3\|=\|Q_3\|_1=0.5H+0.002,\nl
     \disp
     \mu=3\max_{p=\overline{1,3}}\|Q_p\|=3H.
  \end{array}
$$
Next, we calculate
$$
  h_1=\sum_{p=1}^3|\alpha_{p,0}|,\:\:\:
  h_2=\left\{\begin{array}{l}
                \mu h_1^2+(\|A\|+2\mu)h_1,
                \:\:\mbox{if}\:\:h_1>1,\\
                \|A\|+\mu\:\:\mbox{otherwise},
             \end{array}
      \right.
$$
\begin{eqnarray}
   \Delta t=\dfrac{1}{h_2+\delta},\label{dt}
\end{eqnarray}
where $\Delta t$ is the integration step and $\delta$ is an any positive
number (can take a very small).

As seen, the integration step is calculated quite simple
which makes it use in practice. A detailed description of the
algorithm of constructing the approximate chaotic solutions
for the any time segment is given in the next section.

\section{Algorithm for construction of approximate solution}
\label{sec:3}

Before we will seek the approximate solutions of the system
\reff{dynrsys}, it is necessary determine the boundaries what
is limited of researched solution. The sphere $S_a$, limiting the
attractor, may be this boundary. We can set $S_a$, e.g., based on:
\begin{enumerate}
\item On the results of a preliminary calculating experiment by
any numerical method for estimating the approximate boundaries of attractor.

\item The Lyapunov function \cite{Lorenz,Chin} (a classic approach).

\item A theorem on the localization of global attractors \cite{Leonov}.

\item The iteration theorem and first order extremum theorem \cite{Zhang}.

\item A unified approach \cite{Li,PeiWang1,PeiWang2} to estimate
the ultimate bounds of a class of high dimensional quadratic autonomous
dynamical systems.
\end{enumerate}

Let $B_a$ is a ball bounded by the sphere $S_a$.

\begin{figure}[!]
   \centering\includegraphics[width=\textwidth]{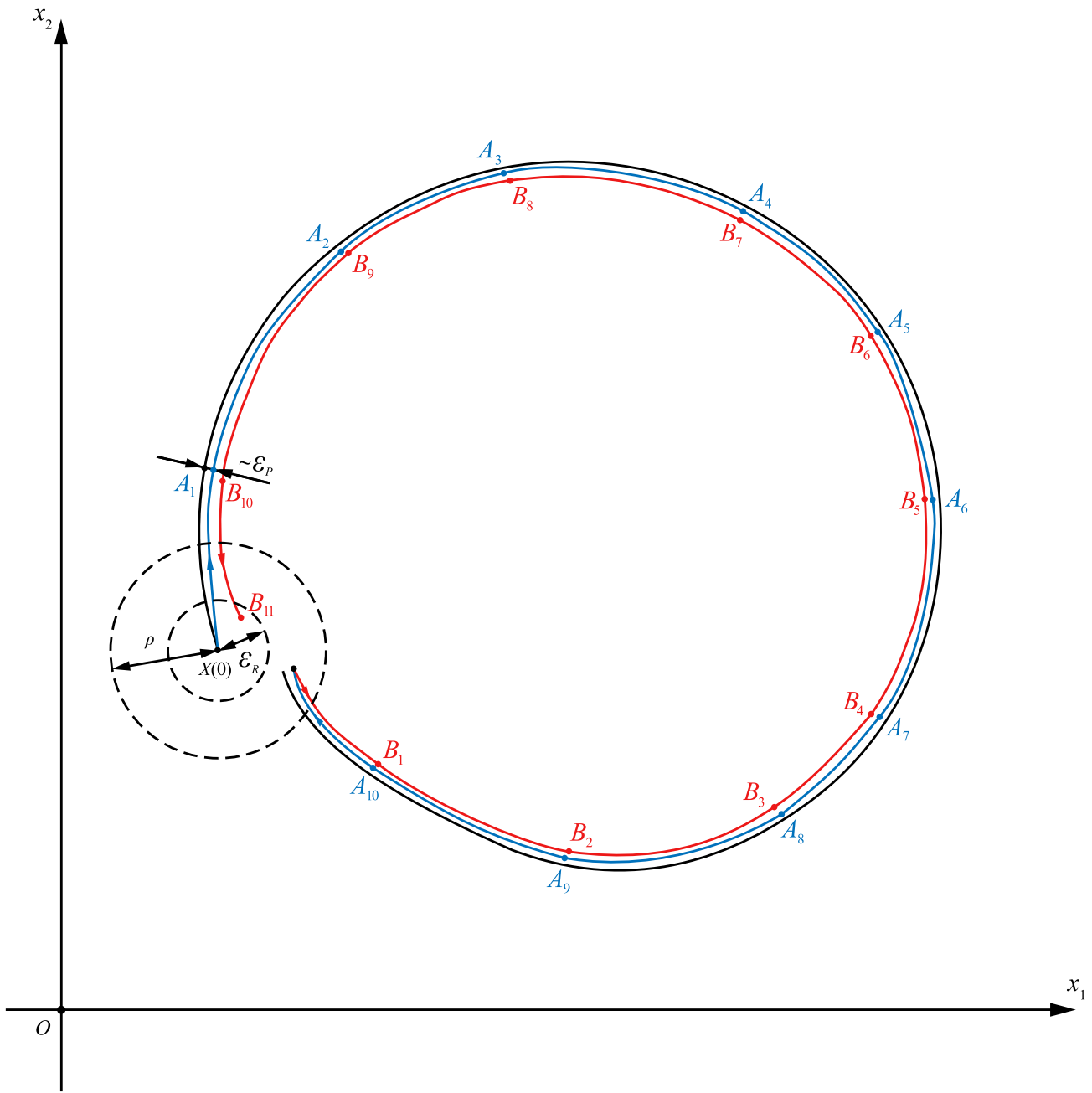}
   \caption{An example of illustration of the FGBFI-method on the plane.
   The points $A_1,\ldots,A_{10}$ connect the parts of the approximate solution
   in forward time where the power series \reff{powser} are convergent,
   the points $B_1,\ldots,B_{10}$ connect the parts for the backward time.
   The black line is projection of the trajectory arc of the
   system \reff{dynrsys} (in general case) on the plane $x_1Ox_2$.
   }\label{fig:1}
\end{figure}

Next, we show the constructing algorithm of approximate solution:
\begin{enumerate}
\item \textbf{Set} the quantity $b_m$ of bits for the mantissa of a real number
and accuracy $\eps_p$ for the power series expansion. The number $b_m$ determines the
machine epsilon $\eps_m$. So we need to take with a reserve of this value, i.e.
choose the value $b_m$, so that
$$
  \eps_m\ll\eps_p;
$$

\item $t:=0$;

\item \textbf{Set} $X(0)\in B_a$ for the system \reff{dynrsys}, $way$ is direction in time
(for going forward $way=1$, going backward $way=-1$), and $T$
(a length of the time segment);

\item $ended:=\mbox{{\ttfamily false}}$;

\item\label{step5} \textbf{Calculate} the integration step $\Delta t$
according to the formula \reff{dt} for $X(0)$;

\item \textbf{If} $\Delta t>T-t$ \textbf{then} $\Delta t:=T-t$, $t:=T$

\textbf{Else} $t:=t+\Delta t$;

\item $\Delta t:=way\cdot\Delta t$;

\item \textbf{Calculate} the point $X(\Delta t)$ with the given accuracy $\eps_p$
for the power series expansions;

\item \textbf{Print}$\big(way\cdot t,\:X(\Delta t)\big)$;

\item \textbf{If} $X(\Delta t)\notin B_a$ \textbf{then}
\textbf{Print}({\sffamily "Decrease the value $\eps_p$ and/or $\eps_m$"}),

$ended:=\mbox{{\ttfamily true}}$;

\item \textbf{If} $t=T$ \textbf{then} $ended:=\mbox{{\ttfamily true}}$;

\item \textbf{If} $ended$ \textbf{then} \textbf{Finish} the algorithm;

\item $X(0):=X(\Delta t)$;

\item \textbf{Goto} Step \ref{step5}.
\end{enumerate}

This algorithm can be applied to forward in time and backward too,
making it a universal.

An example of illustration of the FGBFI-method on the plane is
shown in Fig. \ref{fig:1}. The points $A_1,\ldots,A_{10}$ and
$B_1,\ldots,B_{10}$ are the projections on the plane $x_1Ox_2$
of the points $X(\Delta t)$ when running this algorithm.
If a value of the accuracy $\eps_p$ is large, then
following in the backward time, we will go to infinity,
because the solutions are strongly unstable at $t\rightarrow-\infty$.
Therefore, the algorithm uses the ball $B_a$ for the control
of finding approximate solutions within the boundaries of the attractor.
For return to the given neighborhood $\eps_R$, the value $\eps_p$
(and the number $b_m$ respectively) selected in the calculating experiment.
In fact, the value $\eps_R$ determines how many digits of each coordinate of
the point $B_{11}$ (see Fig. \ref{fig:1}) must coincide with the digits of corresponding
coordinates of the initial point $X(0)$ when we construct the approximate
solution in the backward time. Also we use the configuration analysis of the
approximate chaotic solution to check the accuracy it. In this case, we
calculate the maximum degrees of piecewise polynomials which must be
the same at the forward and backward time as in the articles
\cite{LoziPchelintsev,LozPogPchel}.

\section{Calculating experiments}
\label{sec:4}

We made a calculating experiment for $N=5$, $H=3$ and $I=0.7$
\cite{Mansilla} by the FGBFI-method. In the calculation, the point
$$
  x_1(0)=0.1450756817,\:\:x_2(0)=0.8395885828,\:\:x_3(0)=9.954786333
$$
is found near the attractor. The calculation parameters are
$b_m=160$, $\eps_m=1.36846\cdot10^{-48}$ and $\eps_p=10^{-40}$.
Following by the backward time, it is enough to get the coincidence
of all the decimal places ($\eps_R=10^{-10}$)
of the initial conditions for computing in the time segment $[0,27.327]$.
Also, the maximum degrees of piecewise polynomials coincide
at the forward and backward time, i.e. the criteria
of the article \cite{LozPogPchel} for checking the accuracy of the
approximate solution are performed.

\begin{figure}[!]
   \centering\includegraphics{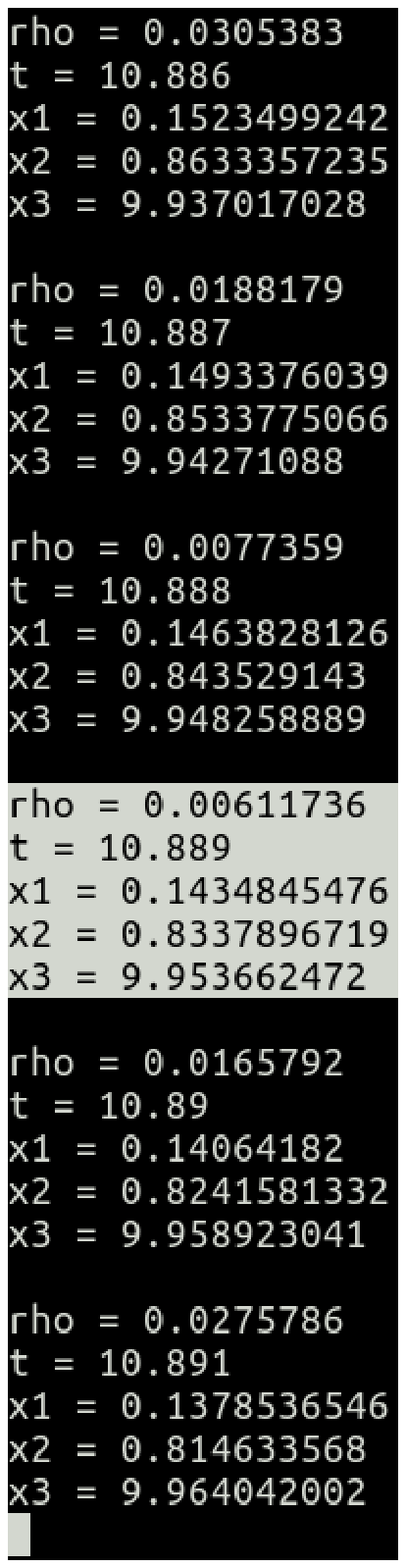}
   \caption{The screenshot of a computer program (for Linux) which
            outputs the distance $\rho$ and coordinates
            $x_1$, $x_2$ and $x_3$ of the trajectory point at
            time $t$ with the time step 0.001.}\label{fig:2}
\end{figure}

\begin{figure}[!]
   \centering\includegraphics[height=8cm]{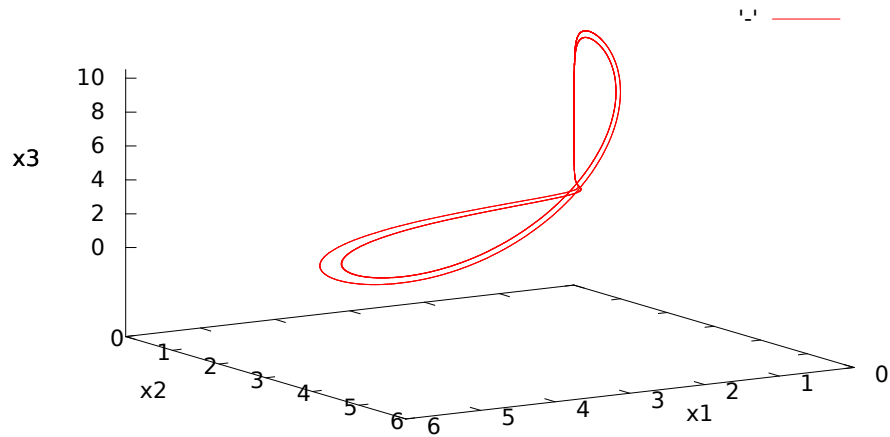}
   \caption{The trajectory arc constructed in the time segment
   $[0,27.327]$ for $N=5$, $H=3$ and $I=0.7$.
   $x_1(0)=0.1450756817$, $x_2(0)=0.8395885828$, $x_3(0)=9.954786333$.}
   \label{fig:3}
\end{figure}

\begin{figure}[!]
   \centering\includegraphics[height=8cm]{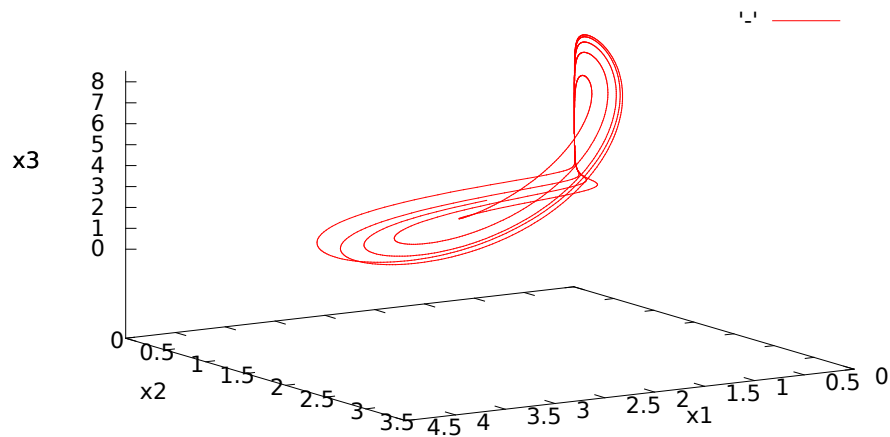}
   \caption{The trajectory arc constructed in the time segment
   $[0,30]$ for $N=5$, $H=3$ and $I=0.4$.
   $x_1(0)=1.292927957$, $x_2(0)=0.5183621413$, $x_3(0)=1.168939477$.}
   \label{fig:4}
\end{figure}

\begin{table}[!]
\caption{The results of calculating experiment.}
\centering\vspace{10pt}
{\begin{tabular}{|c|c|c|c|c|c|}\hline
$n$ & $t_n$ & $x_1(t_n)$ & $x_2(t_n)$ & $x_3(t_n)$ & $\rho(t_n)$\\\hline
0 &      0 & 0.1450756817 & 0.8395885828 &  9.954786333 & 0\\\hline
1 &  5.553 & 0.1201387594 & 0.7151506515 & 9.6198216985 & 0.358201\\\hline
2 & 10.889 & 0.1434845476 & 0.8337896719 &  9.953662472 & 0.006117\\\hline
3 & 16.439 & 0.1207485467 & 0.7178109534 & 9.6243463945 & 0.353004\\\hline
4 & 21.778 & 0.1437352539 & 0.8342333601 & 9.9494643143 & 0.007668\\\hline
5 & 27.327 &  0.118689978 & 0.7111230373 & 9.6323947777 & 0.348049\\\hline
\end{tabular}}
\label{tab:1}
\end{table}

We recorded the rapprochements of trajectory with the initial point
to the Table \ref{tab:1} (with the time step 0.001, see the highlighted strings in
Fig. \ref{fig:2}), wherein $t_0=0$,
$$
  \rho(t)=\sqrt{(x_1(t)-x_1(t_0))^2+(x_2(t)-x_2(t_0))^2+(x_3(t)-x_3(t_0))^2},
$$
since the asymptotic trajectory is Poisson stable. Based on the observed values
$$
  \begin{array}{c}
     \rho(t_0)\approx\rho(t_2)\approx\rho(t_4),\:
     \rho(t_1)\approx\rho(t_3)\approx\rho(t_5),\\
     t_2-t_0\approx t_4-t_2\approx t_3-t_1\approx t_5-t_3\approx 10.89,
  \end{array}
$$
we have an approximation to the periodic solution with the period 10.89.
However, the maximum Lyapunov exponent $\lambda_{\max}$ of this solution
is positive and near to zero, and the Kaplan--Yorke dimension is near to
integer value (see Table \ref{tab:5}). Thus, for large values $t$, we leave
the periodic regime (there is a weak chaotic solution). The trajectory arc
constructed in the time segment $[0,27.327]$ is presented in Fig. \ref{fig:3}.

The chaotic behaviour of the trajectories is observed
$N=5$, $H=3$ and $I=0.4$ (see Fig. \ref{fig:4}). Here, we also got the numerical
solutions of the system \reff{tumorsys} by the 4th order Runge--Kutta (RK4)
method, researching the error $\epsrk$ of this method for different steps
$\dtrk$ (with a constant value) of integration,
$$
  \epsrk=\sqrt{\left(x^{\indrk}_1-x_1(T)\right)^2+
               \left(x^{\indrk}_2-x_2(T)\right)^2+
               \left(x^{\indrk}_3-x_3(T)\right)^2},
$$
where $T=30$, $x^{\indrk}_1$, $x^{\indrk}_2$ and $x^{\indrk}_3$ are
the values of numerical solution by the RK4-method at $t=T$.
The results are shown in Table \ref{tab:2}. Since we are using the 10th characters
after the decimal point as the accurate, then the error
$3.68753\cdot 10^{-9}$ is not so great in relation to the length $T$ of
the integration interval. Here, the maximum degree of the polynomials is
equal to 25, the minimum degree is equal to 15 for the FGBFI-method.

\begin{table}[!]
\caption{The error of the RK4-method compared with the FGBFI-method.}
\centering\vspace{10pt}
{\begin{tabular}{|c|c|c|c|c|c|}\hline
$\dtrk$ & $\epsrk$\\\hline
   0.05 & 0.0387658\\\hline
   0.01 & $4.06488\cdot 10^{-5}$\\\hline
  0.005 & $2.40695\cdot 10^{-6}$\\\hline
  0.001 & $3.68753\cdot 10^{-9}$\\\hline
\end{tabular}}
\label{tab:2}
\end{table}

\begin{table}[!]
\caption{Comparison of the lengths $T$ of the integration intervals and
accuracies $\eps_p$ for different dynamic systems.}
\centering\vspace{10pt}
{\begin{tabular}{|c|c|c|c|c|c|}\hline
{\it Dynamical system} & $T$ & $\eps_p$\\\hline
The Lorenz system \cite{Lorenz,Pchelintsev} & 6.827 & $10^{-50}$\\\hline
The Chen system \cite{Chen,Ueta,LoziPchelintsev} & 8.411 & $10^{-53}$\\\hline
The Sprott--Jafari system \cite{LozPogPchel,Jafari} & 34 & $10^{-15}$\\\hline
The system \reff{tumorsys} \cite{Mansilla} & 30 & $10^{-40}$\\\hline
\end{tabular}}
\label{tab:3}
\end{table}

We compared the lengths $T$ of the integration intervals and accuracies
$\eps_p$ for different dynamic systems in Table \ref{tab:3}. As this table shows,
the value $\eps_p$ is not so small for the Sprott--Jafari system \cite{Jafari}.
This can be explained by the fact that the almost periodic solution (that we
received in the article \cite{LozPogPchel}) has very near to zero (or even negative)
value $\lambda_{\max}$. Also note,  $\lambda_{\max}$ for
solution in Fig. \ref{fig:4} of the system \reff{tumorsys} is positive and near
to zero (see Sect. \ref{sec:5}). Therefore, we have not such a big errors for
the RK4-method.

\section{Calculation of Lyapunov exponents}
\label{sec:5}

Usually, many researchers construct a linearized system of ODEs for
the system \reff{tumorsys} to determine the Lyapunov exponents.
We propose to expand the system \reff{tumorsys} by adding the linearized
equations. The resulting system of 6th-order will also have a quadratic
right-hand side. Let us show it.

Let $x_4(t)$, $x_5(t)$ and $x_6(t)$ are perturbations. We find
(it is assumed that the vector $X$ is made up of three components)
$$
  \begin{array}{c}
  \dfrac{\partial\big(AX+\Phi(X)\big)}{\partial X}
         \left[
            \begin{array}{c}
               x_4\\
               x_5\\
               x_6
            \end{array}
         \right]=\nl
  =
  \left[
     \begin{array}{c}
        2Nx_4-2x_1x_4-Hx_3x_4-Hx_1x_6\nl
        x_1x_4+(4-I)x_5-0.28x_2x_5-0.5Hx_3x_5-0.5Hx_2x_6+0.002x_3x_6\nl
        0.14x_2x_5+0.5Hx_3x_5-Ix_6+0.5Hx_2x_6-0.004x_3x_6
     \end{array}
  \right].
  \end{array}
$$

Now we will work with the extended system \reff{dynrsys}.
Then the matrix $A$ has the form
$$
  A=\left[
       \begin{array}{cccccc}
          2N &   0 &  0 &  0 &   0 &  0\\
           0 & 4-I &  0 &  0 &   0 &  0\\
           0 &   0 & -I &  0 &   0 &  0\\
           0 &   0 &  0 & 2N &   0 &  0\\
           0 &   0 &  0 &  0 & 4-I &  0\\
           0 &   0 &  0 &  0 &   0 & -I
       \end{array}
    \right],
$$
the matrices $Q_1$, $Q_2$ and $Q_3$ will contain zeros
in new places,
$$
  Q_4=\left[
       \begin{array}{cccccc}
           0 & 0 & 0 & -2 & 0 & -H\\
           0 & 0 & 0 &  0 & 0 &  0\\
           0 & 0 & 0 & -H & 0 &  0\\
           0 & 0 & 0 &  0 & 0 &  0\\
           0 & 0 & 0 &  0 & 0 &  0\\
           0 & 0 & 0 &  0 & 0 &  0
       \end{array}
    \right],
$$
$$
  Q_5=\left[
       \begin{array}{cccccc}
           0 & 0 & 0 & 1 &     0 &     0\\
           0 & 0 & 0 & 0 & -0.28 & -0.5H\\
           0 & 0 & 0 & 0 & -0.5H & 0.002\\
           0 & 0 & 0 & 0 &    0  &     0\\
           0 & 0 & 0 & 0 &    0  &     0\\
           0 & 0 & 0 & 0 &    0  &     0
       \end{array}
    \right],
$$
$$
  Q_6=\left[
       \begin{array}{cccccc}
           0 & 0 & 0 & 0 &    0 &      0\\
           0 & 0 & 0 & 0 & 0.14 &   0.5H\\
           0 & 0 & 0 & 0 & 0.5H & -0.004\\
           0 & 0 & 0 & 0 &    0 &      0\\
           0 & 0 & 0 & 0 &    0 &      0\\
           0 & 0 & 0 & 0 &    0 &      0
       \end{array}
    \right],
$$
$$
  \begin{array}{c}
     \|Q_4\|=\|Q_4\|_1=H+2,\nl
     \|Q_5\|=\|Q_5\|_1=\max\{0.5H+0.002,\:1\},\nl
     \|Q_6\|=\|Q_6\|_1=0.5H+0.14,\nl
     \disp
     \mu=6\max_{p=\overline{1,6}}\|Q_p\|=6(H+2).
  \end{array}
$$

Let
$$
  x_4(t)=\sum_{i=0}^\infty\alpha_{4,i} t^i,\:\:\:
  x_5(t)=\sum_{i=0}^\infty\alpha_{5,i} t^i,\:\:\:
  x_6(t)=\sum_{i=0}^\infty\alpha_{6,i} t^i,
$$
$$
  \begin{array}{c}
     \disp
     r_{6,i} =\sum_{j=0}^{i}\alpha_{1,j}\alpha_{4,i-j},\:\:
     r_{7,i} =\sum_{j=0}^{i}\alpha_{3,j}\alpha_{4,i-j},\nl
     \disp
     r_{8,i} =\sum_{j=0}^{i}\alpha_{1,j}\alpha_{6,i-j},\:\:
     r_{9,i} =\sum_{j=0}^{i}\alpha_{2,j}\alpha_{5,i-j},\nl
     \disp
     r_{10,i}=\sum_{j=0}^{i}\alpha_{3,j}\alpha_{5,i-j},\:\:
     r_{11,i}=\sum_{j=0}^{i}\alpha_{2,j}\alpha_{6,i-j},\nl
     \disp
     r_{12,i}=\sum_{j=0}^{i}\alpha_{3,j}\alpha_{6,i-j}.
  \end{array}
$$
Then
$$
  \left\{\begin{array}{l}
     \alpha_{4,i+1}=\dfrac{2N\alpha_{4,i}-2r_{6,i}-Hr_{7,i}-Hr_{8,i}}{i+1},\nl
     \alpha_{5,i+1}=\dfrac{(4-I)\alpha_{5,i}+r_{6,i}-0.28r_{9,i}-
                    0.5Hr_{10,i}-0.5Hr_{11,i}+0.002r_{12,i}}{i+1},\nl
     \alpha_{6,i+1}=\dfrac{-I\alpha_{6,i}+0.14r_{9,i}+0.5Hr_{10,i}+0.5Hr_{11,i}
                           -0.004r_{12,i}}{i+1}.
  \end{array}\right.
$$
We supplement the recurrence relations \reff{recrel} by these relations.

We use the following modification of Benettin's algorithm to determine
the Lyapunov exponents:
\begin{enumerate}
\item \textbf{Divide} the segment $[0,T]$ by segments with length
$\tau=T/M$, $M$ is the quantity that is given;

\item \textbf{Let}
$Y^{(k)}=\left[\alpha_{1,0}^{(k)}\:\alpha_{2,0}^{(k)}\:\alpha_{3,0}^{(k)}\right]$,
$Z^{(k)}_{(1)}=\left[\alpha_{4,0}^{(k,1)}\:\alpha_{5,0}^{(k,1)}\:\alpha_{6,0}^{(k,1)}\right]$,
where $k=\overline{0,M}$. Similarly, \textbf{Introduce} two more vectors $Z^{(k)}_{(2)}$
and $Z^{(k)}_{(3)}$ similarly;

\item \textbf{Input} vector of the initial conditions $Y^{(0)}$ of researched solution
for the system \reff{tumorsys}. \textbf{Input} $Z^{(0)}_{(1)}$, $Z^{(0)}_{(2)}$ and
$Z^{(0)}_{(3)}\neq\textbf{0}$.

\item $k:=0$, $\lambda_1:=0$, $\lambda_2:=0$, $\lambda_3:=0$;

\item\label{step5_} \textbf{If} $k\neq0$ \textbf{then}
$\lambda_1:=\lambda_1+\ln\left|Z^{(k)}_{(1)}\right|$,
$\lambda_2:=\lambda_2+\ln\left|Z^{(k)}_{(2)}\right|$,
$\lambda_3:=\lambda_3+\ln\left|Z^{(k)}_{(3)}\right|$;

\item \textbf{Perform} the normalization
$$
  Z^{(k)}_{(1)}:=\dfrac{Z^{(k)}_{(1)}}{\left|Z^{(k)}_{(1)}\right|};
$$

\begin{table}[!]
\caption{The groups of initial values $Z^{(0)}_{(m)}$
(before normalization) for the linearized system of ODEs.}
\centering\vspace{10pt}
{\begin{tabular}{|c|c|c|c|}\hline
Group number & $Z^{(0)}_{(1)}$ & $Z^{(0)}_{(2)}$ & $Z^{(0)}_{(3)}$\\\hline
I &
$
 \scriptsize\begin{array}{l}
    \alpha_{4,0}^{(0,1)}=5,\\
    \alpha_{5,0}^{(0,1)}=7,\\
    \alpha_{6,0}^{(0,1)}=13
 \end{array}
$
&
$
 \scriptsize\begin{array}{l}
    \alpha_{4,0}^{(0,2)}=10,\\
    \alpha_{5,0}^{(0,2)}=-1,\\
    \alpha_{6,0}^{(0,2)}=11
 \end{array}
$
&
$
 \scriptsize\begin{array}{l}
    \alpha_{4,0}^{(0,3)}=8,\\
    \alpha_{5,0}^{(0,3)}=6,\\
    \alpha_{6,0}^{(0,3)}=9
 \end{array}
$
\\\hline

II &
$
 \scriptsize\begin{array}{l}
    \alpha_{4,0}^{(0,1)}=-6,\\
    \alpha_{5,0}^{(0,1)}=13,\\
    \alpha_{6,0}^{(0,1)}=5
 \end{array}
$
&
$
 \scriptsize\begin{array}{l}
    \alpha_{4,0}^{(0,2)}=63,\\
    \alpha_{5,0}^{(0,2)}=1,\\
    \alpha_{6,0}^{(0,2)}=-17
 \end{array}
$
&
$
 \scriptsize\begin{array}{l}
    \alpha_{4,0}^{(0,3)}=31,\\
    \alpha_{5,0}^{(0,3)}=-7,\\
    \alpha_{6,0}^{(0,3)}=19
 \end{array}
$
\\\hline

III &
$
 \scriptsize\begin{array}{l}
    \alpha_{4,0}^{(0,1)}=1,\\
    \alpha_{5,0}^{(0,1)}=-4,\\
    \alpha_{6,0}^{(0,1)}=75
 \end{array}
$
&
$
 \scriptsize\begin{array}{l}
    \alpha_{4,0}^{(0,2)}=7,\\
    \alpha_{5,0}^{(0,2)}=-13,\\
    \alpha_{6,0}^{(0,2)}=11
 \end{array}
$
&
$
 \scriptsize\begin{array}{l}
    \alpha_{4,0}^{(0,3)}=-40,\\
    \alpha_{5,0}^{(0,3)}=51,\\
    \alpha_{6,0}^{(0,3)}=39
 \end{array}
$
\\\hline

IV &
$
 \scriptsize\begin{array}{l}
    \alpha_{4,0}^{(0,1)}=1,\\
    \alpha_{5,0}^{(0,1)}=1,\\
    \alpha_{6,0}^{(0,1)}=2
 \end{array}
$
&
$
 \scriptsize\begin{array}{l}
    \alpha_{4,0}^{(0,2)}=1,\\
    \alpha_{5,0}^{(0,2)}=-37,\\
    \alpha_{6,0}^{(0,2)}=11
 \end{array}
$
&
$
 \scriptsize\begin{array}{l}
    \alpha_{4,0}^{(0,3)}=29,\\
    \alpha_{5,0}^{(0,3)}=-3,\\
    \alpha_{6,0}^{(0,3)}=5
 \end{array}
$
\\\hline
\end{tabular}}
\label{tab:4}
\end{table}

\begin{table}[!]
\caption{The estimates of Lyapunov exponents and Kaplan--Yorke dimension
for solution in Fig. \ref{fig:3}.}
\centering\vspace{10pt}
{\begin{tabular}{|c|c|c|c|c|}\hline
Group number & $\lambda_1$ & $\lambda_2$ & $\lambda_3$ & $D_{KY}$\\\hline
I   & 0.0233993 & 0.0172255    & $-2.15924$ & 2.0188\\\hline
II  & 0.0433011 & 0.00520866   & $-2.16712$ & 2.0224\\\hline
III & 0.0159841 & $-0.0156199$ & $-2.11898$ & 2.0233\\\hline
IV  & 0.018629  & $-0.0180543$ & $-2.11919$ & 2.0318\\\hline
\end{tabular}}
\label{tab:5}
\end{table}

\begin{table}[!]
\caption{The estimates of Lyapunov exponents and Kaplan--Yorke dimension
for solution in Fig. \ref{fig:4}.}
\centering\vspace{10pt}
{\begin{tabular}{|c|c|c|c|c|}\hline
Group number & $\lambda_1$ & $\lambda_2$ & $\lambda_3$ & $D_{KY}$\\\hline
I   & 0.113902 & $-0.726796$ & $-1.82634$ & 1.1567\\\hline
II  & 0.104372 & $-0.444632$ & $-2.09897$ & 1.2347\\\hline
III & 0.115022 & $-0.472064$ & $-2.08219$ & 1.2437\\\hline
IV  & 0.112198 & $-0.454614$ & $-2.09682$ & 1.2468\\\hline
\end{tabular}}
\label{tab:6}
\end{table}

\item \textbf{Calculate}
$$
  a^{(k)}:=\left\langle Z^{(k)}_{(2)}, Z^{(k)}_{(1)}\right\rangle,\:
  Z^{(k)}_{(2)}:=Z^{(k)}_{(2)}-a^{(k)}Z^{(k)}_{(1)};
$$

\item \textbf{Perform} the normalization
$$
  Z^{(k)}_{(2)}:=\dfrac{Z^{(k)}_{(2)}}{\left|Z^{(k)}_{(2)}\right|};
$$

\item \textbf{Calculate}
$$
  \begin{array}{c}
     b^{(k)}:=\left\langle Z^{(k)}_{(3)}, Z^{(k)}_{(1)}\right\rangle,\:
     c^{(k)}:=\left\langle Z^{(k)}_{(3)}, Z^{(k)}_{(2)}\right\rangle,\nl
     Z^{(k)}_{(3)}:=Z^{(k)}_{(3)}-b^{(k)}Z^{(k)}_{(1)}-c^{(k)}Z^{(k)}_{(2)};
  \end{array}
$$

\item \textbf{Perform} the normalization
$$
  Z^{(k)}_{(3)}:=\dfrac{Z^{(k)}_{(3)}}{\left|Z^{(k)}_{(3)}\right|};
$$

\item \textbf{If} $k\neq M$ \textbf{then} \textbf{Build} the three solutions of
the extended system \reff{dynrsys}
in the time segment $[0,\tau]$ according to the algorithm in
Sect. \ref{sec:3} with forward time. In this case, the initial conditions
$X^{(k)}_{(1)}(0)$, $X^{(k)}_{(2)}(0)$ and $X^{(k)}_{(3)}(0)$
at $k$-th iteration for \reff{dynrsys} formed as
$$
  X^{(k)}_{(m)}(0)=\left[Y^{(k)}\:Z^{(k)}_{(m)}\right]\tr,\:m=\overline{1,3}.
$$
The first three components in the each obtained solution at $t=\tau$
$X^{(k)}_{(m)}(\tau)$ are the same.
\textbf{Record} them in $Y^{(k+1)}$, the other components are recorded in
$Z^{(k+1)}_{(m)}$ respectively;

\item $k:=k+1$;
\item \textbf{If} $k\le M$ \textbf{then Goto} Step \ref{step5_};
\item
$$
  \lambda_1:=\dfrac{\lambda_1}{T},\:
  \lambda_2:=\dfrac{\lambda_2}{T},\:
  \lambda_3:=\dfrac{\lambda_3}{T};
$$
\item \textbf{Print}($\lambda_1$, $\lambda_2$, $\lambda_3$).

\end{enumerate}

We made the computational experiments for the four groups of vectors
$Z^{(0)}_{(m)}$ (see Table \ref{tab:4}). Their results are shown in
Tables \ref{tab:5} and \ref{tab:6} ($D_{KY}$ is the Kaplan--Yorke dimension).
The initial values of the vector components $Z^{(0)}_{(1)}$,
$Z^{(0)}_{(2)}$ and $Z^{(0)}_{(3)}$ (before normalization) for
the linearized system of ODEs are selected randomly.
Also note, $M=20000$. Increasing $M$ does not affect
the given values in Tables \ref{tab:5} and \ref{tab:6}.

A feature of this modification of Benettin's algorithm is to
combine a linearized system of ODEs and
the researched system \reff{tumorsys} in the general form
\reff{dynrsys}.

\section{Acknowledgments}
The reported study was funded by RFBR according to the research project 20-01-00347.

\end{document}